\documentclass[12pt]{article}

\newtheorem{theorem}{Theorem}
\newtheorem{corollary}[theorem]{Corollary}

\newtheorem{lemma}[theorem]{Lemma}
\newtheorem{fact}[theorem]{Fact}

\usepackage{graphicx}
\usepackage{amssymb, amsmath, latexsym}
\begin{document}

\date{}
\title{Distinguishing colorings of Cartesian products of complete graphs}
\author{Michael J. Fisher\thanks{%
Department of Mathematics, California State University, Fresno,
Fresno, CA 93740, email: mfisher@csufresno.edu }
\and Garth Isaak\thanks{%
Department of Mathematics, Lehigh University, Bethlehem, PA 18015,
email: gisaak@lehigh.edu} } \maketitle

\begin{abstract}  We determine the values of $s$ and $t$ for which
there is a coloring of the edges of the complete bipartite graph
$K_{s,t}$ which admits only the identity automorphism. In particular
this allows us to determine the distinguishing number of the
Cartesian product of complete graphs.
\end{abstract}

The distinguishing number of a graph is the minimum number of colors
needed to label the vertices so that the only color preserving
automorphism is the identity. The distinguishing number was
introduced by Albertson and Collins in \cite{ac} and a number of
papers on this topic have been written recently. In this paper we
determine values of $c,s,t$ for which the Cartesian product of
complete graphs of sizes $s$ and $t$ have an identity $c$ coloring.
In particular this allows us to determine the distinguishing number
of the Cartesian product of complete graphs. For $s \leq t$, the
distinguishing number of the Cartesian product of complete graphs on
$s$ and $t$ vertices, $D(K_s \Box K_t)$ is either $\lceil (t+1)^{1/s}
\rceil$ or $\lceil (t+1)^{1/s} \rceil + 1$ and it is the smaller
value for large enough $t$. In almost all cases it can be determined
directly which value holds. In a few remaining cases the value can be
determined by a simple recursion.

Our original version of this paper \cite{fi} was titled `Edge colored
complete bipartite graphs with trivial automorphism groups'. We later
discovered the connection to the distinguishing number.  The current
version has a final added section making the connection to
distinguishing numbers. Thus the rest of paper, except the final
section where we make the connection to distinguishing numbers is the
original version written in terms of identity edge colorings of
complete bipartite graphs.

Harary and Jacobson \cite{hj} examined the minimum number of edges
that need to be oriented so that the resulting mixed graph has the
trivial automorphism group and determined some values of $s$ and
$t$ for which this number exists for the complete bipartite graph
$K_{s,t}$. These are values for which there is a mixed graph
resulting from orienting some of the edges with only the trivial
automorphism. Such an orientation is called an identity
orientation. Harary and Ranjan \cite{hr} determined further bounds
on when $K_{s,t}$ has an identity orientation. They showed that
$K_{s,t}$ does not have an identity orientation for $t \leq
\lfloor \log_3(s-1) \rfloor$ or $t \geq 3^s - \lfloor \log_3(s-1)
\rfloor$ and that it does have an identity orientation for $\lceil
\sqrt{2s} - 3/2 \rceil \leq t \leq 3^s - \lceil \sqrt{2s} - 3/2
\rceil$. In addition they determined exact values when $2 \leq s
\leq 17$. We will show that the first bound is nearly correct.

Observe that a partial orientation of a complete bipartite graph
with parts $X$ and $Y$ has three types of edges: unoriented,
oriented from $X$ to $Y$, and oriented from $Y$ to $X$. We can
more generally think of coloring the edges with some number $c$
of colors. The case $c=1$ is trivial so we will assume throughout
that $c \geq 2$. Automorphisms map vertices to other vertices in
the same part except possibly when $s = t$. So the partial
orientation case corresponds to the case $c = 3$ except possibly
when $s = t$. An identity orientation exists whenever $s=t$ and
we will observe that an identity coloring also exists when $s=t$
except when $s=t = 1$. Thus except for $s=t=1$ our results using
$c = 3$ will correspond to results for the identity orientations
in \cite{hr}. A (color preserving) automorphism is a bijection
from the vertex set to itself with the color of the edge between
two vertices the same as the color of the edge between their
images. An {\em identity coloring} is a coloring with only the
trivial automorphism.

Our main result is
\begin{theorem} \label{t1}
Let $c \geq 2$ and $s \geq 1$ be given integers. When $s \geq 2$
let $x = \lfloor \log_c (s-1) \rfloor$. Then $K_{s,t}$ has an
identity
$c$-edge coloring if and only if exactly one of the following holds: \\
(i) $s = 1$ and $2 \leq t \leq c$  \\
(ii) $2 \leq s \leq c$ and $1 \leq t \leq c^s - 1$ except for
$c = s = t = 2$ \\
(iii) $s > c$ and ${\displaystyle s \leq c^{1+ x} - \lfloor \log_c x
\rfloor} - 2$ and $x+1 \leq t
\leq c^s - x - 1$ \\
(iv) $s > c$ and ${\displaystyle s \geq c^{1+ x} - \lfloor \log_c x
\rfloor}$ and $x+2 \leq t \leq c^s - x - 2$ \\
(v)  $s > c$ and ${\displaystyle s = c^{1+ x} - \lfloor \log_c x
\rfloor} - 1$  and $x+2 \leq t \leq c^s - x - 2$ or $t = x+1,c^s - x
- 1$ and $K_{x+1,s}$ has an identity $c$-edge coloring, except for
the case $c = 2$ and $s = t = 3$.
\end{theorem}

Observe that we can determine if there is an identity coloring
directly from $s$ and $t$ except in case (v) when ${\displaystyle s =
c^{1+ x} - \lfloor \log_c x \rfloor} - 1$ and $t = x+1$ or $t = c^s -
x - 1$. In this situation we let $s' = x + 1$ and $t' = s$ and check
the conditions for $s'$ and $t'$. This needs to be repeated at most
$\log_c^* (s-1)$ times where $\log_c^* (s-1)$ is the iterated
logarithm base $c$.

We are considering the problem of identity orientations examined by
Harary and Ranjan \cite{hr}  in a more general setting and will use a
different notation, however, many of our results  and proofs are
direct extensions of those in \cite{hr}.

A coloring with $c$ colors of the edges of a complete bipartite
graph $K_{s,t}$ having parts $X$ of size $s$ and $Y$ of size $t$
corresponds to a $t$ by $s$ matrix with entries from $\{0,1,
\ldots, c-1\}$. The $i,j$ entry of the matrix is $k$ whenever the
edge between the $i^{th}$ vertex in $Y$ and the $j^{th}$ vertex
in $X$ has color $k$. We will call this the bipartite adjacency
matrix (the usual case being that of general bipartite graph,
which can be thought of as a two coloring, edges and non-edges,
of a complete bipartite graph). For edge colored complete
bipartite graphs, the parts $X$ and $Y$ map to themselves if $|X|
\not = |Y|$. In this case, if $A$ is the bipartite adjacency
matrix, then an automorphism corresponds to selecting permutation
matrices $P_Y$ and $P_X$ such that $A = P_Y A P_X$. If $|X| =
|Y|$ then we also have automorphisms of the form $A = P_Y A^T
P_X$. These will play a role in our results only for certain
small cases. We will discuss these exceptions and then  be able
to consider only  identity colorings  of the form where
 the only solution to $A = P_Y AP_X$ has both $P_Y$ and $P_X$
identity matrices. Throughout we will assume that permutation
matrices are of the correct size for multiplication without stating
the size explicitly. We will use the term identity coloring to refer
to both the edge colored $K_{s,t}$ and to the corresponding adjacency
matrix $A$.

{\bf Results}

We will prove our main theorem  by proving a series of lemmas which
will be stated in terms of the matrix perspective described above.

For a given $c$ and $s$ we will call any $c^s$ by $s$ matrix with
rows corresponding to the $c^s$ distinct $c$-ary $s$-tuples {\em
full}. Any two full matrices of the same size differ only by a
permutation of their rows. If a $c$-ary $t$ by $s$ matrix $A$ has
distinct rows, then its complement ${A^*}$ is `the' $c^s - t$ by
$s$ matrix with rows consisting of the $c$-ary strings of length
$s$ that are not rows of $A$. (The ordering of the rows of
${A^*}$ will not matter for our purposes.)

For any matrix with entries from $\{0,1, \ldots, c-1\}$ the {\em
degree} of a column is a $c$-tuple $(x_0,x_1, \ldots, x_{c-1})$
with $x_i$ equal to the number of entries that are $i$ in the
column. Note that $\sum_{i=0}^{c-1} x_i$ equals the number of
rows of $A$. Thinking of the bipartite adjacency matrix of as
corresponding to a two edge coloring of a complete bipartite
graph then the degree of vertices in $X$ would be $x_1$ in the
degree of $A$. Note that the degree of a vertex and its image in
an automorphism are the same.

The following basic facts, which can easily be checked, will be used.

\begin{fact} \label{f1}
Let $A$ be the adjacency matrix of a $c$-edge colored complete
bipartite graph then: \\
(i) If $A$ is full then so are $P_1A$ and $AP_2$ for permutation
matrices $P_1$, $P_2$ of appropriate sizes. \\
(ii) If there are two identical rows in $A$ then $A$ is not an
identity coloring. \\
(iii) If $A$ is not square and if the columns of $A$ have
distinct degrees and the rows are distinct then $A$ is an
identity coloring. If $A$ is square, has distinct rows, distinct
column degrees and the multiset of column degrees is different
from the multiset of row
degrees then $A$ is an identity coloring. \\
(iv) $A$ is an identity coloring if and only if $A^T$ is.
\end{fact}

We first consider the cases when $|X| = |Y|$.

\begin{lemma} \label{st}
Let $s$ be a non-negative integer. $K_{s,s}$ has an identity
$c$-edge coloring except when $s = 1$ (for any $c$) and when $c =
2$ and $ s= 2$ or $s = 3$.
\end{lemma}

Proof: When $s = 1$ the graph consists of a single colored edge.
Switching the vertices is a non-trivial automorphism.

It is straightforward to check that every 2 coloring of $K_{2,2}$
has a non-trivial automorphism. (Note that when there is exactly
one edge of one of the colors, a non-trivial automorphism must map
$X$ to $Y$ and vice-versa.) It is straightforward, but tedious to
check that every 2 coloring of $K_{3,3}$ has a non-trivial
automorphism.

For $c \geq 3$ the adjacency matrices ${\displaystyle \left[
\begin{array}{cc} 0 & 1 \\ 0 & 2 \end{array} \right] }$ and
${\displaystyle \left[
\begin{array}{ccc} 0 & 1 & 2 \\ 0 & 1 & 0 \\ 0 & 0 & 1 \end{array}
\right] }$ are identity colorings by Fact \ref{f1} (iii).

For $c \geq 2$ and $s \geq 4$ start with a matrix with entries $1$
above the main diagonal and $0$ elsewhere and then replace the 4
by 4 matrix of the first 4 rows and columns with ${\displaystyle
\left[
\begin{array}{cccc} 0 & 1 & 0 & 1 \\ 0 & 1 & 0 & 0 \\ 0 & 0 & 1 &
1 \\ 0 & 0 & 0 & 1
\end{array} \right] }$. This is an identity coloring by Fact \ref{f1}
(iii).  \hfill $\Box$

Note that for identity orientations when $s = t = 1$ the
orientation consisting  of a single edge is an identity
orientation.  This is the only case for complete bipartite graphs
where the existence of an identity orientation and identity
3-edge coloring are not the same.

The next lemma is the $c$ colors version of Lemma 1 of \cite{hj}.

\begin{lemma} \label{l1}
For any integers $c \geq 2$ and $s \geq 2$, $K_{s,t}$ does not have
an identity $c$-edge coloring for $t \geq c^s$.
\end{lemma}

Proof: If at least two rows of the corresponding adjacency matrix are
the same apply Fact \ref{f1} (ii). If not, then $t = c^s$ and the
corresponding adjacency matrix is full. For any non-trivial
permutation matrix $P_2$ (which exist when $s \geq 2$), $AP_2$ is
also full by Fact \ref{f1} (i). Thus, for some permutation matrix
$P_1$ we have $P_1AP_2 = A$. \hfill $\Box$

The $c = 3$ version of the next lemma is used implicitly several
times in \cite{hr}.

\begin{lemma} \label{l2}
Let $A$ be the adjacency matrix of a $c$-edge colored complete
bipartite graph $K_{s,t}$. If $s \not = t$ and $A$ has distinct
rows, then $A$ is an identity coloring if and only if its
complement ${A^*}$ is an identity coloring. In the case that $t =
s$ or $t = c^s - s$, the same holds if we exclude automorphisms
that switch the parts.
\end{lemma}

Proof: Assume that $A$ is not an identity coloring (for $s=t$ assume
that there is a non-trivial automorphism that maps each part to
itself). Since $s \not = t$, for some permutation matrices, $P_1$,
$P_2$, at least one of which is not the identity, we have $A =
P_1AP_2$ and thus $P_1^TA = AP_2$. If $P_2$ is the identity, then
$P_1A = A$ and $P_1$ is also the identity. Thus $P_2$ is not the
identity.

The block matrix
$ {\displaystyle \left[ \begin{array}{c}  A \\
{A^*} \\  \end{array} \right] }$ is full. By Fact \ref{f1}(i)
so is ${\displaystyle  \left[ \begin{array}{c}  A \\ {A^*} \\
 \end{array} \right] P_2 =
\left[ \begin{array}{c}   P_1^T A \\
{A^*} P_2 \\  \end{array} \right] } $. Thus ${A^*}P_2$ and ${A^*}$
have the same rows. Hence for some permutation matrix $P_3$ we
have ${A^*} = P_3 {A^*} P_2$ and thus ${A^*}$ is not an identity
coloring. \hfill $\Box$

\begin{corollary} \label{c1}
For $s \not = t$, $s \geq 2$ and $t \leq c^s$, $K_{s,t}$ has an
identity $c$-edge coloring if and only if $K_{s,c^s-t}$ has an
identity $c$-edge coloring.
\end{corollary}

Proof: Let $A$ be an identity coloring. By Fact \ref{f1} (ii) the
rows are distinct and Lemma \ref{l2} applies. \hfill $\Box$

Applying Lemmas \ref{l1} and \ref{l2} yield the following $c$ color
version of Corollary 7 in \cite{hr}.

\begin{lemma} \label{l3} Let $s \geq 2$.
If $t \leq \lfloor \log_c (s-1) \rfloor$ or $t \geq c^s - \lfloor
\log_c (s-1) \rfloor$ then $K_{s,t}$ does not have an identity
$c$-edge coloring.
\end{lemma}

Proof: The cases $t \geq c^s$ are covered by Lemma \ref{l1}. By
Corollary \ref{c1} it is enough to consider $ t \leq \lfloor
\log_c (s-1) \rfloor$ for the remaining cases. If $A$ is an
identity coloring in these cases, then $A$ has distinct rows and,
by Lemma \ref{l2}, ${A^*}$ is an identity coloring. By Fact
\ref{f1} (iv), the transpose of ${A^*}$ is an identity coloring
with $s' =
 t \leq \lfloor \log_c (s-1) \rfloor$ rows and $t' = s$ columns.
Then $c^{s'} < s = t'$ contradicting Lemma \ref{l1}. \hfill $\Box$

The next lemma is the $c$ color version of Theorem 13 in \cite{hr}
and the proof is similar.

\begin{lemma} \label{l4} Let $s \geq 1$.
Let $r$ be the smallest integer such that ${r + c-1 \choose r} \geq
s$. For $r \leq t \leq c^s - r$ there exists a $c$-ary $t$ by $s$
matrix $A$ with distinct rows and distinct column degrees.
Furthermore, $A$ is an identity coloring except possibly when $A$ is
square.
\end{lemma}

Proof: The furthermore follows immediately from Fact \ref{f1} (iii).

If $A$ has distinct rows and distinct column sums then so does its
complement ${A^*}$. Thus it is enough to prove the theorem for $r
\leq t \leq c^s/2$.

We will use the color set $\{0,1, \ldots, c-1\}$.  For $u \leq v$
let $B_{u,v}$ be the $u$ by $v$ matrix with the $i^{th}$ row
consisting of all zeros except for a $1$ in column $i$ for $1
\leq i < u$ and row $u$ having the first $u-1$ entries $0$ and
the remaining entries $1$. So when $u = 1$ the matrix has one row
of all $1$'s and when $u = v$ it is the identity matrix. Note
that each $B_{u,v}$ and its complement ${B^*_{u,v}}$ has constant
column degrees, has distinct rows and has at least one $1$ in
each column.

Observe that when $s = 1$, $r=0$. For $1 \leq t \leq c$ make the
$i^{th}$ row $i$. (The case $r=0$ can be considered to be true as
there is only one column.) When $s = 2$, $r = 1$. Consider the
$c^2 \times 2$ matrix with rows specified as follows: each $i$ in
$\{1,2, \ldots, c^2\}$ can be written uniquely as $i = jc + k$ for
some $j \in \{0,1, \ldots, c-1\}$ and $k \in \{1,2, \ldots, c\}$.
Let row $i$ be ${\displaystyle \left[ \begin{array}{cc} j & c - k
\end{array} \right] }$. For $1 \leq t \leq c^2 - 1$ taking the
first $t$ rows of this matrix gives the needed matrix as can
easily be checked.

For $s \geq 3$ use induction on $s$. Let $r'$ be the smallest
integer such that ${r' + c - 1 \choose r'} \geq s-1$. Note that
$r' = r$ or $r' = r-1$. For $c = 2$ we have $r = s-1$. For $c \geq
3$, by the choice of $r$ we have $s > {(r-1) + c -1 \choose c-1}
\geq {r-1 + 2 \choose 2}$ which implies that $r < \sqrt{2s}$.

For $r \leq t \leq c^s/2$ we will consider several cases. Note that
$c^s/2 \leq (c-1)c^{s-1}$ for $c \geq 2$. For $c = 2$ cases 1 and 2
suffice.

Case 1: $r \leq t \leq c^{s-1} - r$. Since $r' \leq r$, by
induction there exists a $t$ by $s-1$ matrix with distinct rows
and distinct column degrees. The $s-1$ column degrees each satisfy
$x_0 + x_1 + \cdots + x_{c-1} = t$ where $t \geq r$. So there
exists a solution $x^*_0 + x^*_1 + \cdots + x^*_{c-1} = t$
distinct from any of the degrees. Add a new column $s$ with
$x^*_i$ entries equal to $i$. The rows are still distinct and the
new column degree is also distinct from the first $s-1$.

Case 2: $(a+1)c^{s-1} - r < t \leq (a+1)c^{s-1}$ for some
non-negative integer $a$. Let $u = t - ((a+1)c^{s-1} - r)$ so $1 \leq
u \leq r$.

For $c = 2$: Then $a=0$ since we need only consider $t \leq
2^s/2$.
 The cases $s=2$ and $s=3$ are easily checked.
Consider $s > 3$. Let $D$ be the $s-1$ by $s-1$ matrix with zeros
above the main diagonal and ones elsewhere. Then $D^*$ has
distinct rows, distinct column degrees and at least $2^{s-2} -
(s-2) \geq 2$ zeros in each column. Take the rows of $B_{u,s-1}$
with a last column of $0$'s added along with the rows of $ D^*$
with a last column of $1$'s added. The result has distinct rows.
The last column has $u$ zeros and the other columns have at least
$2$ zeros from $D^*$ and $u-1$ zeros from $B_{u,s-1}$ for a total
of more than $u$. Thus the last column degree is distinct from the
others. Since also the column degrees of $ D^*$ are distinct and
those of $B_{u,s-1}$ are constant we get distinct column degrees.

For $c > 2$: Let $A'$ be a solution with $s-1$ columns and $t =
c^{s-1} - r$ rows, which exists by induction. Note that $a \leq c-2$
since if $a \geq c-1$ then $t > ((c-1)+1)c^{s-1} - r \geq c^s -
\sqrt{2s} \geq c^s/2$.

If $a=0$ take $A'$ with a last column of $0$'s added and
$B_{u,s-1}$ with a last column of $2$'s added. Each part has
distinct rows and the last entries for the rows differ for the
different parts so the rows are distinct. On the first $s-1$
columns $B_{u,s-1}$ has constant degree and $A'$  has distinct
degrees so the first $s-1$ column degrees are distinct. The last
column has no $1$'s and every other column does since $B_{u,s-1}$
does. Thus the last column has degree distinct from the others.

Now assume that $a \geq 1$. Let $D$ be the $ac^{s-1}$ by $s$
matrix with rows all $c$-ary $s$-tuples with last entry from
$\{2,3, \ldots, a+1 \}$. $D$ has constant degree on the first
$s-1$ columns and distinct rows. Take $D$, $B_{u,s-1}$ with a last
column of $0$'s added and $A'$ with a last column of $1$'s added.
Each part has distinct rows and the last entries for the rows
differ for the different parts so the rows are distinct. On the
first $s-1$ columns $D$ and $B_{u,s-1}$ with the appended column
have constant degree and $A'$ with the appended column has
distinct degrees so the first $s-1$ column degrees are distinct.
The number of $0$'s on each of the first $s-1$ columns is at least
$(u-1) + ac^{s-2}$ from the $0$'s in $B_{u,s-1}$ and $D$
respectively. This is strictly greater than $u$ as $c \geq 3$, $s
\geq 3$ and $a \geq 1$. The last column has $u$ $0$'s. Thus the
last column has degree distinct from the others.

Case 3: $(a+1)c^{s-1} <  t \leq (a+1)c^{s-1} + r$ for some
non-negative integer $a$. From the remarks above we can also
assume that $c \geq 3$. Let $u = t - (a+1)c^{s-1}$ so $1 < u \leq
r$. Note that $a \leq c - 3$ since if $a \geq c-2$ then $t >
(c-1)c^{s-1} \geq c^s/2$. Since $r < \sqrt{2s}$ one can then check
that $u+r \leq 2r \leq c^{s-1} - r$. Thus, by induction there is a
solution $A'$ with $s-1$ columns and $u+r$ rows. Let $D$ be the
$ac^{s-1}$ by $s$ matrix with rows all $c$-ary $s$-tuples with
last entry from $\{2,3, \ldots, a+1 \}$. In the case that $a = 0$,
$D$ will be empty. $D$ has constant degree on the first $s-1$
columns and distinct rows. Take $D$, $B_{r,s-1}$ with a last
column of $0$'s added and $A'$ with a last column of $(c-1)$'s
added. Since $a \leq c-3$, color $(c-1)$ is not used on $D$.
Hence, each part has distinct rows and the last entries for the
rows differ for the different parts so the rows are distinct. The
last column has no $1$'s and every other column does since
$B_{u,s-1}$ does. Thus the last column has degree distinct from
the others.

Case 4: $(a+1)c^{s-1} + r < t < (a+2)c^{s-1} - r$ for some
non-negative integer $a$. Let $u = t - (a+1)c^{s-1}$ so that $r  < u
< c^{s-1} - r$. Note that $a \leq c-3$ since if $a \geq c-2$ then $t
> (c-1)c^{s-1} + r \geq c^s - c^{s-1} \geq c^s/2$.
Let $A'$ be a solution with $s-1$ columns and $u$ rows which
exists by induction. Let $D$ be the $(a+1)c^{s-1}$ by $s$ matrix
with rows all $c$-ary $s$-tuples with last entry from $\{2,3,
\ldots, a+2\}$. $D$ has constant degree on the first $s-1$ columns
and distinct rows and has $c^{s-1}$ $1$'s in each of these
columns. Since $s \geq 2$ and $c \geq$ there is at least one $1$
in each of the first $s-1$ columns. Take $D$ and $A'$ with a last
column of $0$'s added. Each part has distinct rows and the last
entries for the rows differ for the different parts so the rows
are distinct. On the first $s-1$ columns $D$ has constant degree
and $A'$  has distinct degrees so the first $s-1$ column degrees
are distinct. The last column has no $1$'s and every other column
does since $D$ does on these columns. Thus the last column has
degree distinct from the others.
 \hfill $\Box$

Observe that this lemma is best possible in the sense for $t < r$
no such set could exist. If it did  then we would have $s > {t+c-1
\choose t}$ distinct solutions to $x_0 + x_1 + \cdots + x_{c-1} =
t$, a contradiction. By complementation no such set exists for $t
> c^s - r + 1$.

{\bf Proof of Theorem \ref{t1}}:

For (i) in the theorem, the case $s=t=1$ is noted in lemma \ref{st}.
If $t > c$ then two edges have the same color and by Fact \ref{f1}
(ii) the coloring is not an identity coloring. For $2 \leq t \leq c$
assigning different colors to the edges gives an identity coloring by
Fact \ref{f1} (iii).

For the remaining cases we use induction on $s$.

For $s \leq c$ we have $r = 1$ in Lemma \ref{l4} and (ii) follows
except when $s = t$. The $s = t$ cases are covered by Lemma \ref{st}.
The $s \leq c$ cases will also be the basis for the induction.

Lemma  \ref{l3} covers the cases $t \geq c^s - x$ and $t \leq x$.
Lemma \ref{l4} covers the cases $s \leq t \leq c^s - s$ since $r \leq
s$ always. For the remaining cases it is enough to consider $x+1 \leq
t < s$ by Corollary \ref{c1}.

In each case we will let $s' = t$ and $t' = s$ and use fact \ref{f1}
(iv), that $K_{s,t}$ has an identity coloring if $K_{s',t'}$ does.
Then since $s' = t < s$ we can inductively check $K_{s',t'}$.

For $x+2 \leq t < s$ we need to show that $K_{s',t'}$ has an identity
coloring. Note that $c^{s'} \geq c^{x+2} = (c)c^{1 + \lfloor \log_c
(s-1) \rfloor} \geq c(s-1)$. When $s' \leq c$ we have $c^{s'} - 1
\geq c(s-1)-1 \geq s = t'$ and we get an identity coloring. For $s'
> c$ we have  ${\displaystyle c^{s'} - \lfloor \log_c (s'-1) \rfloor
- 2 \geq  c(s-1) - \lfloor \log_c (s'-1) \rfloor - 2 \geq s = t'}$
and again get an identity coloring.

The case $s' = t = x+1$ remains. When $s = c^{x+1} - \lfloor \log_c x
\rfloor - 1$ statement (v) is that we look at $K_{s',t'}$. We need to
show that $K_{s',t'}$ does not have an identity coloring when (iv) $s
\geq c^{x+1} - \lfloor \log_c x \rfloor$ and it does when (iii) $s
\leq c^{x+1} - \lfloor \log_c x \rfloor - 2$.

First note that if $s' = x+1 \leq c$ then $\lfloor \log_c x \rfloor
\leq \lfloor \log_c (c-1) \rfloor = 0$. So (iv) only occurs if $s =
c^{1+x}$. Then $t' = s = c^{s'} > c^{s'} - 1$ and by (ii) $K_{s',t'}$
does not have an identity coloring. (The case $s' = t = 1$ is already
covered by (i).) When $s' = x+1 \leq c$ and (iii) occurs we have $t'
= s \leq c^{s'} - 2$ and we have an identity coloring.

Now assume that $s' = x+1 > c$. If (iv) then $t' = s \geq c^{x+1} -
\lfloor \log_c x \rfloor = c^{s'} - \lfloor \log_c (s'-1) \rfloor$
and hence by Lemma \ref{l3} $K_{s',t'}$ does not have an identity
coloring. If (iii) then $t' = s \leq c^{x+1} - \lfloor \log_c x
\rfloor - 2 = c^{s'} - \lfloor \log_c (s'-1) \rfloor - 2$ and hence
by induction $K_{s',t'}$ has an identity coloring. \hfill $\Box$

While our Theorem gives an exact answer for determining if $K_{s,t}$
has an identity coloring in nearly all cases (a recursive check is
required in (v)) we give here a few specific examples of using
Corollary \ref{c1} directly to determine if there is an identity
coloring for illustration. We will take $c = 3$. The conclusions in
\cite{hr} state that $K_{s,t}$ has an identity coloring if and only
if $1 \leq t \leq 3^s - 1$ when $s \in \{2,3\}$, if and only if $2
\leq t \leq 3^s - 2$ for $s \in \{4,5, \ldots, 8\}$, and if and only
if $3 \leq t \leq 3^s - 3$ for $s \in \{9,10, \ldots, 17 \}$. Since
$K_{3,t}$ has an identity coloring for $t \in \{1,2, \ldots, 26\}$ we
get, for example, that $K_{26,3}$ has an identity coloring. In a
similar manner we can conclude that $K_{s,t}$ has an identity
coloring if and only if  $3 \leq t \leq 3^s - 3$ for $s \in \{9,10,
\ldots, 26\}$. The facts that $K_{4,79}$ has an identity coloring and
that $K_{3,79}$ does not have an identity coloring and other similar
cases show us that $K_{79,t}$ has an identity coloring if and only if
$4 \leq t \leq 3^{79} - 4$.

{\bf Distinguishing Numbers}

Distinguishing numbers of Cartesian products have been investigated
in \cite{al}, \cite{ik}, \cite{kz}. We have also recently discovered
\cite{ijk} which appears to contain results on the distinguishing
number similar to our paper although our paper seems to have a
slightly narrower range of values for when a recursion is needed to
determine the distinguishing number.

Recall that line graphs of a complete bipartite graphs are Cartesian
products of complete graphs. That is, $L(K_{s,t}) = K_s \Box K_t$.
Thus our results correspond to vertex colorings of $K_s \Box K_t$.
Our automorphisms are on the vertex set of the bipartite graphs
$K_{s,t}$ so we need to observe that they do indeed correspond to
automorphisms of $K_s \Box K_t$. This follows directly from the
following result of Imrich and Miller, cited in \cite{al}. Theorem:
If $G$ is connected and $G = H_1 \Box H_2 \Box \cdots \Box H_r$ is
its prime decomposition, then every automorphism of $G$ is generated
by the automorphisms of the factors and the transpositions of
isomorphic factors. As the cases $s=t$ are easily dealt with we can
directly translate our results on coloring $L(K_{s,t})$ to coloring
$K_s \Box K_t$.

To determine the distinguishing number of $K_s \Box K_t$ we need to
determine the smallest $c$ in Theorem \ref{t1} for which $K_{s,t}$
has an identity coloring. When $s=1$  we have $K_1 \Box K_t = K_t$
and we see from part (i) that $c = t$. This corresponds to the known
result that $D(K_t) = t$. From the upper bounds on $t$ we see that
$c$ should be approximately $\lceil (t+1)^{1/s} \rceil$. If $c <
\lceil (t+1)^{1/s} \rceil$ then $t > c^s - 1$ and Theorem \ref{t1}
tells us that there is no identity coloring. If $c = \lceil
(t+1)^{1/s} \rceil + 1$ then $t+1 \leq (c-1)^s$ and hence $t \leq c^s
- sc^{s-1} - 1 \leq c^s - x - 2$ for $x = \lfloor \log_c(s-1)
\rfloor$ and Theorem \ref{t1} tells us that there is an identity
coloring. So the distinguishing number is $\lceil (t+1)^{1/s} \rceil$
or $\lceil (t+1)^{1/s} \rceil + 1$. In particular, for large $t$
relative to $s$, (for example $t \geq s^s$) we get that the
distinguishing number is $\lceil (t+1)^{1/s} \rceil$.

With the details of Theorem \ref{t1} we get the following. Observe
that in all but one case we determine the distinguishing number
immediately. In the remaining case we determine it from a recursion
 that is similar to that of Theorem \ref{t1}. In particular it is repeated
only when we need to determine if $D(K_{x+1,s}) = c$. So $c$ stays
fixed in the recursive computations and thus number of steps in the
recursion is at most iterated logarithm $\log_c^* (s-1)$.

\begin{corollary}
For $2 \leq s \leq t$ let $c = \lceil (t+1)^{1/s} \rceil$. So $t \leq
c^s - 1$.  Then
 $D(K_s \Box K_t)$ equals $c$ or $c+1$. When $t \geq s^s$ the value
 is $c$.

 In particular, letting
 $x = \lfloor \log_c (s-1) \rfloor$ we have: \\
 (i) $D(K_s \Box K_t) = c$ for $s \leq t \leq c^s - x - 2$ except
 for the case $D(K_2 \Box K_2) = 3$. \\
 (ii) $D(K_s \Box K_t) = c + 1$ for $c^s - x  \leq t \leq c^s - 1$.
 \\
 (iii) $D(K_s \Box K_t) = c$ for $t = c^s - x - 1$ and
 ${\displaystyle s \leq c^{1+ x} - \lfloor \log_c x
\rfloor} - 2$. \\
 (iv) $D(K_s \Box K_t) = c + 1$ for $t = c^s - x - 1$ and
 ${\displaystyle s \geq c^{1+ x} - \lfloor \log_c x
\rfloor}$. \\
(v) When $t = c^s - x - 1$ and
 ${\displaystyle s = c^{1+ x} - \lfloor \log_c x
\rfloor} - 1$ then $D(K_s \Box K_t) = c$ if $D(K_{x+1} \Box K_s) \leq
c$ and $D(K_s \Box K_t) = c + 1$ if $D(K_{x+1} \Box K_s) \geq c + 1$

\end{corollary}

{\em Acknowledgements:} The authors would like to thank Peter Hammer
for encouraging the writing of this paper. Garth Isaak would like to
thank the Reidler Foundation for partial support of this research.

\end{document}